\documentclass [11pt]{article}
\usepackage{amsfonts}
\usepackage{amsmath}
\newtheorem{theorem}{Theorem}[section]
\newtheorem{proposition}{Proposition}
\newtheorem{lemma}{Lemma}
\newtheorem{corollary}{Corollary}
\begin{document}
\newcommand{\TL}{\tilde{\mathcal{L}^{\mathcal{I}}}}
\newcommand{\LL}{\mathcal{L}^{\mathcal{I}}}
\newcommand{\KK}{\mathcal{K}_{\mathcal{I}}}
\newcommand{\iast}{\ast_{\mathcal{I}}}
\newcommand{\es}{\left( \begin{array}{c} g\ast\phi_s \\ \phi_s\end{array} \right)}
\newcommand{\ei}{\left(\begin{array}{c} g\ast\phi_i \\ \phi_i\end{array}
\right)}
\newcommand{\ees}{\left(\begin{array}{c} g\ast_{\mathcal{I}}\phi_s \\ 0 \end{array}
\right)}
\newcommand{\eei}{\left(\begin{array}{c} g\ast_{\mathcal{I}}\phi_i \\ 0 \end{array}
\right)}
\newcommand{\eees}{\left(\begin{array}{c} 0 \\ \phi_s\end{array}
\right)}
\newcommand{\eeei}{\left(\begin{array}{c} 0 \\ \phi_i\end{array}
\right)}
\topmargin 0pt
\headsep 0pt
\date{}
\title{On the Largest Singular Values of Random Matrices with Independent Cauchy Entries}

\author{ 
Alexander Soshnikov\thanks{
Department of Mathematics,
University of California at Davis, 
One Shields Ave., Davis, CA 95616, USA.
Email address: soshniko@math.ucdavis.edu.
Research was supported in part by the Sloan Research Fellowship and the 
NSF grant DMS-0103948. }
and
Yan V. Fyodorov\thanks{
Department of Mathematical Sciences,
Brunel University, Uxbridge UB83PH, UK.
Email address: yan.fyodorov@brunel.ac.uk.}}

\date{}
\maketitle
\begin{abstract}
We apply the method of determinants to study the distribution of the largest singular values
of large $ m \times n $ real rectangular random matrices with  independent  Cauchy entries.
We show that statistical properties of the (rescaled by a factor $ \frac{1}{m^2\*n^2} $) 
largest singular values agree in the limit 
with the statistics of the Poisson random point process with the intensity $ \frac{1}{\pi} \* x^{-3/2} $ and, therefore, 
are different from the Tracy-Widom law.
Among other corollaries of our method  we show an interesting connection between the mathematical expectations
of the determinants  of the complex rectangular $ m \times n $ standard Wishart ensemble and the
real rectangular $ 2m \times 2n $ standard Wishart ensemble.
\end{abstract}

\section{Introduction and Formulation of Results.}

The main goal of this paper is to study the spectral properties of a 
large random matrix with i.i.d. Cauchy 
entries. In other words we consider a rectangular 
$ m\times n $ matrix $ A=(a_{ij}), \ 1 \leq i \leq m, \ 
1 \leq j \leq n, $ where $ \{ a_{ij} \} $ are independent identically distributed Cauchy random variables
with the probability density $ f(x)= \frac{1}{\pi (1+x^2)}$. 
Our goal is to study the singular values of $A$
as the dimensions of a matrix go to infinity, $ m \to \infty, \ n \to \infty.$
This is clearly equivalent to studying the eigenvalues of a positive-
definite $ n \times n $ matrix $ M = A^t\* A.$ Matrices of such type are quite often
called sample covariance 
matrices in random matrix literature.
Positive-definite matrices are of particular importance 
in statistics (we refer to \cite{Mui}, \cite{Wil}, 
\cite{Ja} for the classical works on
statistical applications of the spectral properties of Wishart
matrices, and \cite{Janik} for a few recent developments and
applications to various fields). They also are of long-standing interest 
in nuclear physics, starting with the classical works 
\cite{Wig3}, \cite{Br}. More recently they were used to model the 
"dissipative" part of the effective Hamiltonian in 
quantum chaotic scattering (see \cite{FS} and references therein)
and appeared to be very intimately connected 
with the "chiral" ensembles studied 
in Quantum Chromodynamics, see \cite{qcd}. 
As other important applications of random positive-definite 
matrices we mention that they
are used in a branch of condensed matter theory known as 
mesoscopics to model famous universal conductance fluctuations and other 
transport properties of small metallic 
samples and quantum dots, see \cite{Been}, and also emerged in
theory of information communication in random environment \cite{zeitouni}.
 
It is well known that if the entries of $A$ are i.i.d. random variables 
with zero mean and finite variance $ \sigma^2 $, the empirical distribution function of the eigenvalues of 
$ \frac{1}{n} \* A^t\* A $ converges in the limit $ m \to \infty, \ n \to \infty, \ m/n \to \gamma \in 
(0, +\infty) $ to the Marchenko-Pastur law (see e.g. \cite{MP}, \cite{Bai}) defined by its density
\begin{equation}
\label{PasturMarchenko}
p_{\gamma}(x)= (2\pi x \gamma \sigma^2)^{-1} 
\* \sqrt{(b-x)(x-a)}, \ \ \  a\leq x \leq b ,
\end{equation}
where $ a= \sigma^2(1-\gamma^{-1/2})^2 $ and $ b= \sigma^2(1+\gamma^{-1/2})^2 $ (we assume here 
$ \gamma \geq 1 $). Since the spectrum of $ A^t \* A $ differs from the spectrum of $A \* A^t $ only by the multiplicity
of the
eigenvalue $ \lambda=0$ (for $ m \geq n $ the matrix $ A \* A^t $ has $ m-n $ additional zero eigenvalues)
for the rest of the paper we can assume $ m \geq n .$ 
 Under the assumption that
the fourth moment of $a_{ij}$ is finite Yin , Bai and Krishnaiah
(\cite{YBK}, see also \cite{Sil}, \cite{BaS}) showed that the largest eigenvalue of $ \frac{1}{n} \* A^t\* A $
converges to $b$ almost surely. Recently Johnstone proved that in the standard Wishart case (i.e. $ \{ a_{ij} \}
$ are i.i.d. $ N(0,1) $ random variables) the properly rescaled largest eigenvalue converges in distribution to
the $\beta=1$ Tracy-Widom distribution ((\cite{TW2}), see also (\cite{TW1}). Soshnikov (\cite{So2}) 
generalized the result of Johnstone to the non-Gaussian case  provided $ n-m= O(n^{1/3}) $ and the moments of 
the matrix entries do not grow  very fast.
There are quite a few standard methods that have been successfully used for Wigner and sample
covariance matrices
in the case when second and higher moments of matrix entries exist, most notably the method of
moments (\cite{Wig1}, \cite{Wig2},
\cite{SiS1}, \cite{SiS2}, \cite{So1}, \cite{So2}),
the method of resolvents
(\cite{MP}, \cite{KKP}, \cite{Bai}), the method of orthogonal polynomials (\cite{D}, \cite{J}), the method used by Johansson  \cite{Jo1} 
(and recently extended by Ben Arous and P\'{e}ch\'{e} to the case of sample covariance matrices)  which is based on the 
Kazakov-Br\'{e}zin-Hikami trick (\cite{Kaz},\cite{BH1}, \cite{BH2}), etc.
Unfortunately, the above mentioned approaches are not suitable for 
the Cauchy case. In particular, one can expect the spectral properties of
$A^t\* A $ in the Cauchy case to be rather
different from the case of a finite variance.  In our view this makes the studies of the Cauchy case 
especially interesting. Denote the eigenvalues of $A^tA$ by $ \lambda_1 \geq \lambda_2 \geq \lambda_3 \geq
\ldots \geq \lambda_n. \ $ 
It is expected that the majority of the eigenvalues
are proportional to $ m\*n. $ We would like to specifically single out the reference \cite{CB} where Cizeau and Bouchard studied the 
spectral properties of the Wigner random matrices with the heavy tails
(see also references \cite{Bur}, \cite{Janik1} for physical papers on the so-called L\'{e}vy-Smirnov unitary ensembles). 
Among other things Cizeau and Bouchard argued (on a physical level of rigor) that the empirical distribution function of the 
eigenvalues of a properly normalized Wigner matrix  (with the heavy tails of the marginal distribution of matrix entries)
converges to a limiting distribution that can be obtained as a solution of a quite complicated system of two integral equations 
(we refer to the formulas (15), (12a), (12b) in \cite{CB}). It is not difficult to guess the right order of the normalization: 
after the normalization, the norm of any given matrix raw has to be  of order of constant (in particular in the Cauchy case one has to 
normalize the matrix entries by
$n^{-1}$  and  in the case of a finite variance the 
normalization 
is $ \sqrt{n}$).
The support of the limiting distribution 
is the whole real line. One can expect (arguing at the same level of physical rigor) to derive
a similar system of integral equations for the limiting distribution 
function of the eigenvalues of $  \frac{1}{m\*n} \* A^t \* A $ in the case of the i.i.d.
Cauchy entries of $A$. The support of the limiting distribution should be the positive half of the real line.  Our results formulated below 
indicated that the asymptotics of the spectral density at infinity should be $ \frac{1}{\pi \* x^{3/2}}. $ 
This suggests that the largest
eigenvalues of $ A^t\*A$ grow faster than $m\*n$, in fact we will show below that the largest eigenvalues are of the order of $m^2\*n^2 $. 
Let us rescale the eigenvalues by that factor:

\begin{equation}
\label{rescale}
\tilde{\lambda_i}= \frac{\lambda_i}{n^2 \* m^2}, \ \ i=1, \ldots, n.
\end{equation}
The goal of this paper is to study the local distribution of the largest eigenvalues by the method of determinants.
Our main results are Theorems 1.1 and 1.2 formulated below.

\begin{theorem}
 Let  $A $ be a random rectangular $m \times n$ matrix ($ m \geq n $ ) with i.i.d. Cauchy entries and
$ z $ a complex number with a positive real part.
Then as $ n \to \infty $
we have
\begin{eqnarray}
\label{glavform}
\lim_{n \to \infty} E \left(\det(1 +\frac{z}{m^2 \*n^2} \* A^t \* A )\right)^{-1/2} & =&
E \prod_{i=1}^n (1+ z\* \tilde{\lambda_i})^{-1/2} =
\exp\left(-\frac{2}{\pi} \* \sqrt{z} \right)\\
\label{glavforma}
&= & {\bf E} \prod_{i=1}^{\infty} (1+ z\* x_i)^{-1/2},
\end{eqnarray}
where
we consider the branch of $\sqrt{z}$ on $D=\{ z: \Re z >0 \} $ such that $ \sqrt{1}=1,$
$ E $ denotes the mathematical expectation with respect to the random matrix ensemble defined above, ${\bf E}$ denotes the mathematical 
expectation with respect to the inhomogeneous Poisson random point process on the positive half-axis with the intensity 
$ \frac{1}{\pi \* x^{3/2}}, $ 
and the convergence 
is uniform inside $D$ (i.e. it is unform on  the compact subsets of $D$). For  a real positive $z=t^2, \ t \in \Bbb R, $ one can  estimate 
the rate of convergence, namely
\begin{equation}
\label{glavformA}
\lim_{n \to \infty} E \left(\det(1 +\frac{t^2}{m^2 \*n^2} \* A^t \* A )\right)^{-1/2} =
E \prod_{i=1}^n (1+ t^2\* \tilde{\lambda_i})^{-1/2} =
\exp\left(-\frac{2}{\pi} \* |t| \* \big(1+o(n^{-1/2+\epsilon})\big)\right),
\end{equation}
where $ \epsilon $ is an arbitrary small positive number and the convergence is uniform on the compact subsets of
$ [0, +\infty).$

\end{theorem}
 We will discuss the properties of Poisson random point processes at the Appendix. A very useful introduction to the 
elementary theory and methods of random point processes is \cite{DVJ}. It is not a coincidence that the intensity of the Poisson random 
process in the above theorem is equal to the leading term of the asymptotics of the density of the square of a standard Cauchy random 
variable.

 We claim that the result can be generalized to the case of 
a sparse random matrix with Cauchy entries.  Let as before $ \{ a_{jk}\}, \ 1 \leq j \leq m, \ 1 \leq k \leq n, \ $ be i.i.d. Cauchy random 
variables, and $ Q= (q_{jk}) $ be a $ m \times n $ non-random rectangular $ 0-1$ matrix such that the number of non-zero entries in each 
column
is fixed and equals to $b_n.$ For technical reasons we assume that $b_n $ grows to infinity as some power of $ n$, i.e. 
$ \  \ b_n \geq n^{\alpha}, $ for some $ 0<\alpha \leq 1 $  and $ \ln m $ is much smaller than any power of $n $.
We define a $ m \times n $ rectangular matrix $A$ with the entries
$ \Gamma_{jk}= q_{jk} \* a_{jk}, \ \ 1\leq j \leq m, \ 1 \leq k \leq n. \ $ As before we denote by 
$ \lambda_1 \geq \lambda_2 \ldots \geq \lambda_n $
the eigenvalues of $ \Gamma^t\* \Gamma $. The appropriate rescaling for the largest eigenvalues in this case is going to be
$ \tilde{\lambda_i} = \frac{\lambda_i}{m^2 \* b_n^2}, \ \ i=1, \ldots, n.$
We claim that the result of the Theorem 1.1 can be extended to the case of a sparse random matrix $\Gamma $. 

\begin{theorem}
 Let  $\Gamma $ be a sparse random rectangular $m \times n$ matrix ($ m \geq n $ ) defined as above and
$ z $ a complex number with a positive real part.
Then as $ n \to \infty $
we have
\begin{eqnarray}
\label{glavform1}
\lim_{n \to \infty} E \left(\det(1 +\frac{z}{m^2 \* b_n^2} \* 
\Gamma^t \* \Gamma )\right)^{-1/2} &=&
E \prod_{i=1}^n (1+ z\* \tilde{\lambda_i})^{-1/2} =
\exp\left(-\frac{2}{\pi} \* \sqrt{z}  \right) \\
&= & {\bf E} \prod_{i=1}^{\infty} (1+ z\* x_i)^{-1/2},
\end{eqnarray}
where, as in Theorem 1.1,
we consider the branch of $\sqrt{z}$ on $D=\{ z: \Re z >0 \} $ such that $ \sqrt{1}=1,$
$ E $ denotes the mathematical expectation
with respect to the random matrix ensemble defined in the paragraph above the theorem, ${\bf E}$ denotes the mathematical 
expectation with respect to the inhomogeneous Poisson random point process on the positive half-axis with the intensity 
$ \frac{1}{\pi \* x^{3/2}},$
and the convergence 
is uniform inside $D$ (i.e. it is unform on the compact subsets of $D$). 
For  a real positive $z=t^2, \ t \in \Bbb R, $ one can get an estimate 
on the rate of convergence, namely
\begin{equation}
\label{glavform1A}
E \left(\det(1 +\frac{t^2}{m^2 \* b_n^2} \* \Gamma^t \* \Gamma )\right)^{-1/2} =
E \prod_{i=1}^n (1+ t^2\* \tilde{\lambda_i})^{-1/2} =
\exp\left(-\frac{2}{\pi} \* t \* \big(1+o(b_n^{-1/2+\epsilon})\big)\right),
\end{equation}
where
$ \epsilon $ is an arbitrary small positive number and the convergence is 
uniform on the compact subsets of $[0, +\infty).$

\end{theorem}
The result of Theorem 1.2 can be generalized even further.  Let the setting be as in Theorem 1.2 but relax the condition
that the number of non-zero entries in each column is exactly $b_n$  to the condition $ \sum_{k=1}^n q_{jk} = b_n\* (1+o(1)), \ \ 
j=1, \ldots, m,$  (for example the relaxed condition is satisfied by a typical realization of a random matrix $Q$ with independent Bernoulli
$ 0-1$ entries with $ Pr (q_{jk}=1)= b_n/n ).\ $ Then  we still have
$ E \left(\det(1 +\frac{z}{m^2 \* b_n^2} \* \Gamma^t \* \Gamma )\right)^{-1/2} =
\exp\left(-\frac{2}{\pi} \* \sqrt{z} \right) \* (1+o(1)). \ $ 
The proof is almost identical to the proof of Theorem 1.2 and will be left to the 
reader.

The case when the number of non-zero terms of $Q$ is fixed in 
each raw can be treated in a similar manner. 

An important consequence of Theorems 1.1 and 1.2  is that the statistical properties of the largest eigenvalues $ \tilde{\lambda_1}, 
\tilde{\lambda_2},\ldots $ are drastically different from the statistical properties of the (rescaled) largest eigenvalues
in GOE (\cite{TW2}) and real Wishart case (\cite{J}) that are described by the ($\beta=1$) Tracy-Widom law.

Theorems 1.1 and 1.2  follow from the Proposition 1  formulated below.

\begin{proposition}
Let $A=(a_{jk}) $ be a random rectangular $m \times n$ matrix with independent (not necessarily identically distributed) 
entries with the 
characteristic functions of the matrix entries $ g_{jk}(s)= E \exp(is\*a_{jk}).$ 
Let $ t_i>0, \ i=1, \ldots, r $ be some positive parameters.
Then the following formula holds
\begin{eqnarray}
\label{proddet}
& & E \left( \prod_{i=1}^r \det(1+ t_i^2\* A^t\*A) \right)^{-1/2} = (2\pi)^{-r(n+m)/2} \times \nonumber \\
& &  \int_{R^{r(n+m)}} 
\prod_{i=1}^r d^ns^{(i)}\*
d^mp^{(i)} \ \exp\left( -\sum_{i=1}^r (|s^{(i)}|^2/2 +|p^{(i)}|^2/2) \right) 
\prod_{1\leq j \leq m, 1 \leq k \leq n}
g_{jk}\left(\sum_{i=1}^r t_i \* p_j^{(i)} \* s_k^{(i)}  \right).
\end{eqnarray}
\end{proposition}

Theorems 1.1  and 1.2  imply several important corollaries given in the next section.

The proofs of Theorems 1.1 and 1.2 and Proposition 1 are given  in the next section.
Section 3 is devoted to  application of the method of determinants to random matrices with i.i.d. 
complex entries. We prove in section 3 the analogue of Proposition 1 (Proposition 2) in the complex case and as a corollary 
establish an interesting connection between the determinants in the $ 2m \times 2n $ rectangular real Wishart 
case and $ m \times n $ rectangular complex Wishart case (see Lemma 1).

\section{Proofs of Theorems 1.1 and 1.2}

We start with the proof of the Proposition 1.
Consider 
$ \bigl(\det(1 + t^2 \* A^tA)\bigr)^{-1/2}
, \ \ t>0.$
Let  $ s= (s_1, \ldots, s_n)^t, \ \ p=(p_1, \ldots, p_m)^t \ $ be real $n-$ and $m-$ 
dimensional column vectors.  Let 
$ B(t)= \left( \begin{array}{cc} Id & t i \* A\\ t i \* A^t & Id \end{array} \right) \ $ and
$ \ d^{n}s \* d^{m}p = \prod_{i=1}^n d s_i  \* \prod_{j=1}^m   p_i . \ $

Then
\begin{eqnarray}
\label{propder}
& &\bigl(\det(1 + t^2\* A^tA)\bigr)^{-1/2} = 
\left( \det \left( \begin{array}{cc} 1 &  ti \* A\\ t i \* A^t & 1 \end{array} \right) \right)^{-1/2} 
\nonumber\\
& & = \left(\frac{1}{\pi}\right)^{(n+m)/2} \* \int \ d^{n} s \  d^{m} p  \ 
\exp \bigl(-(s, p) \* B(t) \* (s,p)^t \bigr)   
\nonumber\\
& & = \left(\frac{1}{\pi}\right)^{(n+m)/2} \* 
\int \ d^{n} s \  d^{m} p \  \exp( - (|s|^2 + |p|^2)) \* E \left( \exp (
-2i \sum_{1 \leq j \leq m, 1 \leq k \leq n} t \* a_{jk} \* p_j\* s_k  ) \right).
\end{eqnarray}

In (\ref{propder}) we used standard 
properties of  the Gaussian integral.
The formula (\ref{propder}) and the independence of matrix entries imply

\begin{eqnarray}
\label{propder1}
& & E \left( \prod_{i=1}^r \det(1+ t_i^2\* A^t\*A) \right)^{-1/2} =  \nonumber\\
& &  \pi^{-r(n+m)/2} \int_{R^{r(n+m)}} 
\prod_{i=1}^r d^ns^{(i)}\*
d^mp^{(i)} \ \exp\left( -\sum_{i=1}^r (|s^{(i)}|^2 +|p^{(i)}|^2) \right) \times \nonumber \\ 
& & E \prod_{1\leq j \leq m, 1 \leq k \leq n}
\exp\left(2i \* a_{jk} \* \sum_{i=1}^r t_i \* p_j^{(i)}\*s_k^{(i)}  \right) = (2\pi)^{-r(n+m)/2} \times
\nonumber \\
& &  \int_{R^{r(n+m)}} 
\prod_{i=1}^r d^ns^{(i)}\*
d^mp^{(i)} \ \exp\left( -\sum_{i=1}^r (|s^{(i)}|^2/2 +|p^{(i)}|^2/2) \right) 
\prod_{1\leq j \leq m, 1 \leq k \leq n}
g_{jk}\left(\sum_{i=1}^r t_i \* p_j^{(i)} \* s_k^{(i)} \right).
\end{eqnarray}
The Proposition is proven.

To prove Theorem 1.1 we observe that the functions at the l.h.s. of 
(\ref{glavform}) are analytic and uniformly bounded in $D= \{z: \Re z >0\}
$. Therefore, by the Vitali's theorem it is enough to prove the convergence for real positive $z$.  Let us denote $ z=t^2,$ where $t$ is a 
positive real number and apply the result of the proposition in 
the case $r=1$. Since the matrix entries of $A$ are i.i.d. Cauchy we have
$ g_{jk}(s)= g(s)=\exp(-|s|) $ and 
\begin{eqnarray*}
& & E \left(  \det(1+ t^2\* A^t\*A) \right)^{-1/2} =  \nonumber\\
& &  (2\pi)^{-(n+m)/2} \int_{R^{n+m}} 
d^n s\*
d^m p \ \exp\left( - (|s|^2 +|p|^2)/2 \right) \times 
\prod_{1\leq j \leq n, 1 \leq k \leq m}
\exp\left( -t \*|s_k \* p_j| \right) =\nonumber \\ 
& &  2^m \* \int_{R^n} d^n s \* (2 \pi)^{-n/2} \exp\left( -\frac{1}{2} |s|^2\right)
\* \int_{R^m_{+}} (2 \pi)^{-m/2} \* \exp\left( -\frac{1}{2} \* \sum_{j=1}^m (p_j^2 + 2 p_j \* t \* \sum_{k=1}^n
|s_k|) \right)
\end{eqnarray*}
\begin{equation}
\label{nol}
= \int_{R^n} d^n s \* (2 \pi)^{-n/2} \exp\left( -\frac{1}{2} |s|^2\right)
\* \Psi^m\bigl(t \* \sum_{k=1}^n |s_k|\bigr),
\end{equation}
where $ \Psi(y) = 2 \* e^{y^2/2} \* \int_y^{+\infty} \frac{1}{\sqrt{2\pi}} \* e^{-t^2/2} \* dt .$ In particular,
$\Psi(0)=1 $ and $\Psi'(0)= - (2/\pi)^{1/2}.$  It is easy to see that the function $ \Psi(y) $ is monotonly decreasing on $ [0, +\infty)$, 
in particular $ 1=\Psi(0)=\max_{[0,+\infty)} \Psi(y).$ Indeed, $ \Psi'(y)= 2 \* y \* e^{y^2/2} \* 
\int_y^{+\infty} \frac{1}{\sqrt{2\pi}} \* e^{-t^2/2} \* dt - 2 \* \frac{1}{\sqrt{2\pi}}. $ The assertion then follows from the inequality
$ y^{-1} \* e^{-y^2/2}  > \int_y^{+\infty}  e^{-t^2/2} \* dt $ for $ y>0.$

Replacing $t$ by $\frac{t}{n\*m}$ we arrive at
\begin{equation}
\label{odin}
E \left(  \det(1+ \frac{t^2}{n^2 \* m^2} \* A^t\*A ) \right)^{-1/2}= \int_{R^n} d^n s \* (2 \pi)^{-n/2} 
\exp\left( -\frac{1}{2} |s|^2\right) \* \Psi\Big(\frac{t}{n\*m} \* \sum_{k=1}^n |s_k|\Big)^m.
\end{equation}
The r.h.s. of the last formula suggests to use the law of large numbers and large deviations estimates for
the sum of the  absolute values of $n $ standard Gaussian random variables.
Since $ \int \frac{1}{2\pi} |s| e^{-s^2/2} ds = (2/\pi)^{1/2},$
we see that  $ \frac{1}{n} \* \sum_{j=1}^n |s_j| = (2/\pi)^{1/2} + o(n^{-1/2+\epsilon}) $ with probability 
$ 1 - O(\exp(- n^{1.99 \* \epsilon}).$
Recalling that $\Psi(0)=1 $ and $ \Psi'(0)= - (2/\pi)^{1/2}$ we get
\begin{eqnarray}
\label{dva}
& & E \left(  \det(1+ \frac{t^2}{n^2\*m^2}\* A^t\*A ) \right)^{-1/2}= \int_{R^n} d^n s \* (2 \pi)^{-n/2} 
\exp\left( -\frac{1}{2} |s|^2\right) \* \left(1 - \frac{2t \* (1+o(n^{-1/2 + \epsilon}))}{\pi\* m}\right)^m =
\nonumber \\
& & \left(1 - \frac{2t \* (1+o(n^{-1/2 +\epsilon})}{\pi\* m}\right)^m = 
\exp\left(-\frac{2 \* t \*  (1+o(n^{-1/2 +\epsilon})}{\pi} \right),
\end{eqnarray}
for any $ \epsilon >0 .$
Theorem 1.1 is then follows by the Vitali theorem.

The proof of Theorem 1.2 is very similar. Again we can restrict our attention to the case when $z$ is a real positive number,
$z=t^2.$ We have
\begin{eqnarray}
\label{nolnol}
& & E \left(  \det(1+ t^2\* \Gamma^t\*\Gamma) \right)^{-1/2} =  \nonumber\\
& & \int_{R^n} d^n s \* (2 \pi)^{-n/2} \exp\left( -\frac{1}{2} |s|^2\right)
\* \prod_{j=1}^m \Psi\bigl(t \* \sum_{k=1}^n q_{jk} \* |s_k|\bigr).
\end{eqnarray}

Let $ k_1^{(j)}, k^{(j)}_2, \ldots, k^{(j)}_{b_n} $ be the indices $k $ for which $ q_{jk}=1. $ Then 
$ \sum_{k=1}^n q_{jk} \* |s_k| = \sum_{l=1}^{b_n} |s_{k^{(j)}_l}| $ and we can claim that
$ \frac{1}{b_n} \* \sum_{l=1}^{b_n} |s_{k^{(j)}_l}| = (2/\pi)^{1/2} + o(b_n^{-1/2+\epsilon}) $ with probability 
$ 1 - O(\exp(- b_n^{1.99 \* \epsilon}) \ $ for each $ 1\leq j \leq m $.
Since we assumed that $ \ln(m) $ is much smaller than any power of $n$ we get
similarly to (\ref{dva}) that
\begin{eqnarray}
\label{dva1}
& & E \left(  \det(1+ \frac{t^2}{m^2\*b_n^2}\* \Gamma^t\* \Gamma) \right)^{-1/2}= \int_{R^n} d^n s \* (2 \pi)^{-n/2} 
\exp\left( -\frac{1}{2} |s|^2\right) \* \left(1 - \frac{2t \* (1+o(b_n^{-1/2 + \epsilon}))}{\pi\* m}\right)^m =
\nonumber \\
& & \left(1 - \frac{2t \* (1+o(b_n^{-1/2 +\epsilon})}{\pi\* m}\right)^m = 
\exp\left(-\frac{2 \* t \*  (1+o(b_n^{-1/2 +\epsilon})}{\pi} \right).
\end{eqnarray}
Theorem 1.2 is proven.

Below we restrict our attention to the corollaries of Theorem 1.1 (full matrix case).  In corollaries of Theorem 1.2 are basically identical
to those of Theorem 1.1 (with an obvious change of $n$ to $b_n$ where it is needed).

{\bf Remark 0}
We are not aware that the results of Theorem 1.1 and 1.2 are enough to imply that ths statistics of the largest eigenvalues are Poisson in 
the limit of  $n \to \infty$.  Indeed, to prove the Poisson statistics in the limit one has to show that
\begin{equation}
\label{rav}
\lim_{n \to \infty} E \prod_{i=1}^n \bigl(1+f(\tilde{\lambda_i})\bigr) = {\bf E} \prod_{i=1}^{+\infty}\bigl(1+f(x_i)\bigr) 
\end{equation}
for a sufficiently large class of the test functions $f$. The results of Theorems 1.1-1.2 claim that
(\ref{rav}) is valid for $f(x)= (1+z\*x)^{-1/2} -1 $ for all $z$ such that $ \Re z >0.$ Below we formulate several corollaries of our main 
result that are weaker than the claim about the Poisson statistics, but still give us some information about the behavior of the largest 
eigenvalues. The proof of the Poisson statistics for the largest eigenvalues in Wigner random matrices with heavy tails will appear in 
\cite{So3}.

{\bf Remark 1}

It follows immediately from the result of the Theorem 1.1  that ``only a finite number'' of the eigenvalues
$ \lambda_i $ are of the order of $ n^2 \* m^2.$  Indeed, let $N_{n,m} $ be an integer growing to infinity arbitrary slow as 
$ n \to \infty  $  and $\delta >0 $ be an arbitraty small positive number. Then $ Pr \Big( \#( \lambda_i \geq \delta \* n^2 \* m^2 ) 
\geq N_{n,m}
\Big) \to 0 $ as $ n\to \infty.$ Indeed, suppose this is not the case. Then
$ Pr \Big( \#( \lambda_i \geq \delta \* n^2 \* m^2 )\geq N_{n,m} \Big) \geq \kappa >0 $ and 
$ E \left(  \det(1+ \frac{t^2}{n^2 \* m^2}\* A^t\*A ) \right)^{-1/2} \leq \kappa \* (1+ t^2 \* \kappa)^{-N_{n,m}} + (1-\kappa)\* 1.$
One obtains a contradiction since for $ t= N^{-1}_{n,m} $ the r.h.s. of the last inequality does not go to zero (see also next remark).
One can also rewrite the statement of this remark in the following way: for any positive $ \delta, \ \kappa $  there exist
$n_0(\delta, \kappa)$ and $ C(\delta, \kappa) $ such that
$ \Pr \Big( \#( \lambda_i \geq \delta \* n^2 \* m^2 ) 
\geq C \Big) < \kappa $ for all $n_0 \leq n \leq m .$

{\bf Remark 2}

It is clear from the proof of the Theorem 1.1 that the asymptotic result  
$$
E \left(  \det(1+ \frac{t^2}{n^2 \* m^2}\* A^t\*A ) \right)^{-1/2}=
\exp\left(-\frac{2 \* t \*  (1+o(n^{-1/2 +\epsilon}))}{\pi} \right)
$$
holds uniformly in $t$ on compact subsets
of $ [0, +\infty).$ In particular the result is valid for a sequence $ t_n \to 0.$  

Below we formulate and prove some additional consequences of Theorem 1.1. Our first observation is that one can repeatedly 
differentiate  (\ref{glavform}) with respect to parameter $ z.$

\begin{corollary}
Let $ \Re z >0 .$ Then
\begin{eqnarray}
\label{chetyre}
& & \lim_{n \to \infty} E \prod_{i=1}^n (1+ z\* \tilde{\lambda_i})^{-1/2}
\* \left( \sum_{j=1}^n \frac{\tilde{\lambda_i}}{1+ z\* \tilde{\lambda_i}} \right)=
\frac{2}{\pi} \*  z^{-1/2} \* \exp\left(-\frac{2}{\pi} \*  \sqrt{z} \right) = \\
& & {\bf E} \prod_{i=1}^{\infty} (1+ z\* x_i)^{-1/2}
\* \left( \sum_{j=1}^{\infty} \frac{x_i}{1+ z\* x_i} \right),
\end{eqnarray}
where, as above, {\bf E} stands for the mathematical expectation with respect to the inhomogeneous Poisson random point process on 
$ (0, +\infty) $ with the intensity $\rho(x)= \frac{1}{\pi \* x^{3/2}}.$
\end{corollary}

{\bf Remark 3}

If we let $ z \to 0 $ in (\ref{chetyre}) one gets
$ E \left( \sum_{j=1}^n \tilde{\lambda_i} \right)= +\infty, $
which trivially follows from the fact that matrix entries of $A$ are Cauchy random variables. Essentially the result of the corollary can 
tell
us how fast the mathematical expectation $ E \left( \sum_{j} \tilde{\lambda_j} \right) $ grows if we restrict the summation only to 
$\tilde{\lambda_i} \leq L $ where $  L $ is large.

{\bf Proof of the Corollary 1.}
The result immediately follows form the uniform convergence of the analytic functions in (\ref{glavform}-\ref{glavforma}).

By differentiating (\ref{glavform}) twice one can obtain in a similar fashion that
\begin{eqnarray}
\label{chetyrepljus}
& & \lim_{n\to \infty} E \prod_{i=1}^n (1+ z\* \tilde{\lambda_i})^{-1/2}
\* \left( \left( \sum_{j=1}^n \frac{\tilde{\lambda_i}}{1+ z\* \tilde{\lambda_i}} \right)^2 +
\sum_{i=1}^n \frac{(\tilde{\lambda_i})^2}{(1+z\*\tilde{\lambda_i})^2} \right)= \nonumber \\
& & \left( \frac{4}{\pi^2} \*  z^{-1} + \frac{2}{\pi} \*  z^{-3/2} \right) \* 
\exp\left(-\frac{2}{\pi} \* \sqrt{z} \right).
\end{eqnarray}

\begin{corollary}
 There is a contant $C$ which depends on $ \gamma $ such that for 
$ Pr (\frac{\lambda_1}{n^2\*m^2}  > x ) < C \* x^{- 1/2}  $ uniformly for large 
$ n \leq m $ and $ x $.
\end{corollary}

Indeed, it follows from Theorem 1.1 and Remark 1 that 
$ E (1 + t^2 \* \tilde{\lambda_1})^{-1/2}  
\geq \exp\big(-2\* t \* (1 +o(n^{-1/2+\epsilon}))/\pi\big) $ uniformly 
in $ t $ on compact subsets of $ [0, \infty) $.  Therefore
$ (1 - Pr (\tilde{\lambda_1} >x)) + \frac{1}{\sqrt{1+t^2\* x}}\* Pr(\tilde{\lambda_1} >x) \geq \exp\Big(-\frac{2}{\pi}\* t
\*(1+o(n^{-1/2+\epsilon}))
\Big),$
which implies
$ 1 - \exp\Big(-\frac{2}{\pi}\*  t\*(1+o(n^{-1/2+\epsilon})) \Big)\geq 
\left(1- \frac{1}{\sqrt{1+t^2\* x}}\right) \* 
Pr(\tilde{\lambda_1} >x) .$ Choosing $ t^2 \* x =1 $ we obtain that
$ \frac{2}{\pi} \* \gamma \* x^{-1/2} \* (1+\delta)  
\geq (1- 2^{-1/2}) \* Pr(\tilde{\lambda_1} >x) $ for all 
sufficiently large $n, \ m, $ and $x$.

{\bf Remark 4}

It is not difficult to show that in probability $ \lambda_1 = O(n^2\*m^2). $ 
To see this we observe that the operator norm 
$ \| A \| $ can be bounded from below by $ \max_{1\leq j \leq m, 1\leq k \leq n } |a_{jk}| .$ The maximum of 
$n\*m $ i.i.d. Cauchy random variables is of the order $ O(n\*m) $ (with the limiting distribution of
$ \frac{1}{n\*m} \*\max_{1\leq j \leq m, 1\leq k \leq n } |a_{jk}| $ easily computable, namely
$ \Pr \left(\frac{1}{n\*m} \*\max_{1\leq j \leq m, 1\leq k \leq n } |a_{jk}|  \leq x \right) \to \exp(-\frac{2}{\pi\*x})$).
We expect that the 
limiting distribution $ \frac{\lambda_1}{n^2 \*m^2}  $  also exists and is given by the distribution of the rightmost particle in the 
Poisson process with the intensity $\rho(x)= \frac{1}{\pi} \* x^{-3/2},$ in other words,
$\lim_{n \to \infty} \Pr (\tilde{\lambda_1} <x) = \exp( -\frac{2}{\pi} \* x^{-1/2}).$

It is a useful excercise to see what Proposition 1 gives in the
Wishart case. Below we treat the case of one 
determinant.
\begin{eqnarray*}
& & E \left(  \det(1+ t^2\* A^t\*A) \right)^{-1/2} =  \\
& &  (2\pi)^{-(n+m)/2} \int_{R^{n+m}} 
d^n s\*
d^m p \ \exp\left( - (|s|^2 +|p|^2)/2 \right) \times 
\prod_{1\leq j \leq n, 1 \leq k \leq m}
\exp\left( -\frac{1}{2} (t \*s_k \* p_j)^2 \right) =  \\
& &  \* \int_{R^n} d^n s \* (2 \pi)^{-n/2} \exp\left( -\frac{1}{2} |s|^2\right)
\* \int_{R^m} d^m p \* (2 \pi)^{-m/2} \* \prod_{k=1}^m \exp\left( -\frac{1}{2} \* p_k^2 \* \Big(1+  t^2 \* \sum_{j=1}^n
s_j^2 \Big) \right)= \\
& & \int_{R^n} d^n s \* (2 \pi)^{-n/2} \exp\left( -\frac{1}{2} |s|^2\right)
\* \big(1+ t^2 \* \sum_{j=1}^n s_j^2 \big)^{-m/2} 
\end{eqnarray*}
\begin{eqnarray}
\label{gaussnol}
&= & c_n \* \int_0^{+\infty} \exp(-\frac{1}{2} r^2) \* r^{n-1} \* (1+t^2\*r^2)^{-m/2} \*dr  \nonumber \\
&=& c_n \* 2^{\frac{n}{2}-1} \* \int_0^{+\infty} 
e^{-r} \* r^{\frac{n}{2} -1} \* (1+2\* t^2\*r)^{-m/2} dr,
\end{eqnarray}
where $c_n$ is the normalization constant, $ c_n^{-1} =\int_0^{+\infty} \exp(-\frac{1}{2} r^2) \* r^{n-1} dr=
2^{\frac{n}{2}-1} \* \Gamma(n/2).$ To study the global distribution of the eigenvalues in the Wishart ensemble
one has to consider rescaling $t^2 \to  \frac{t^2}{n}$ (since typical eigenvalues of $A^t\*A$ are of the order of $n$). It follows 
that
\begin{eqnarray}
\label{gaussodin}
& & E \left(  \det(1+ \frac{t^2}{ n} \* A^t\*A) \right)^{-1/2} =  
\left[\Gamma(n/2)\right]^{-1} \* \int_0^{+\infty} e^{-r} \* r^{\frac{n}{2} -1} \* (1+ 2 \* t^2\*\frac{r}{n})^{-m/2} \*dr 
= \nonumber \\
& & \left[\Gamma(n/2)\right]^{-1} \* (n/2)^{n/2} \int_0^{+\infty} 
e^{-\frac{n}{2} \* r} \* r^{\frac{n}{2}} 
\* (1+ t^2\*r)^{-m/2} \*r^{-1} \*dr= \nonumber \\
& & \left[\Gamma(n/2)\right]^{-1} \* (n/2)^{n/2} \* \int_0^{+\infty} e^{-\frac{n}{2}\*\mathcal{L}(z)} \* z^{-1} \* dz, 
\end{eqnarray}
where $ \ \mathcal{L}(z)=  z + \frac{m}{n}\* \ln (1+  t^2 \* z ) -  \ln z. \ $
The asymptotics of the last integral can be obtained by the steepest descent method. The formulas are especially
simple
in the square case
$ m=n $. One then can find a positive solution of the equation
\begin{equation}
\frac{d\mathcal{L}}{dz}= 1+ \frac{t^2}{1+ t^2 z} - \frac{1}{z} = 
\frac{t^2 \* z^2 + z - 1}{ z\*(1 + t^2 \* z)}=0
\end{equation}
to be $ \ z(t)= \frac{-1 + \sqrt{4 \*t^2 +1}}{2 \* t^2}. \ $ 
Taking into account that
$ \ \frac{d^2\mathcal{L}}{dz^2} = \sqrt{4 \* t^2+1} \ $ at $ \ z=z(t) \ $ 
we obtain that in the square Wishart case
\begin{eqnarray}
\label{predel}
& &  E \left( \det(1 + \frac{t^2}{n} \* A^t \* A) \right)^{-1/2} = \nonumber \\
& & \frac{2^n \* n^{n/2}}{ \Gamma(n/2)} \* 
\exp \left( -n \* 
\frac{-1 + \sqrt{4 \*t^2 +1}}{4 \*t^2}
 \right)
\* 
(1+\sqrt{4\*t^2+1})^{-n}
\* 
\sqrt{\frac{2 \pi \* \sqrt{4 \* t^2+1}}{n }} \* 
\frac{2 t^2}
{-1 + \sqrt{4 \* t^2 +1}}
 \* (1 +o(1)).
\end{eqnarray}
The fact that the asymptotics in (\ref{predel}) is exponential in $n$ is  standard . Indeed, it is a  straightforward excercise to 
verify that
\begin{equation}
\label{zbc}
\lim_{n \to \infty} \frac{1}{n}\* \ln \left( E \left(\det(1 + \frac{t^2}{n} \* A^t \* A) \right)^{-1/2}\right) = 
-\frac{1}{2} \* \int_0^4 (\ln(1 + t^2\*x) \* p_1(x)\* dx ,
\end{equation}
where $ p_1(x) $ is the probability density of the Marchenko-Pastur law defined in (\ref{PasturMarchenko}), which reflects
the law of large numbers for the linear statistics $ \sum_{i=1}^n  \ln(1 + t^2\* \lambda_i\*n^{-1})  $
where $ \{\lambda_i\}_{i=1}^n $ are the eigenvalues of the real Wishart matrix. The variance of the linear statistics is bounded and has a 
limit as $ n \to \infty $ (for a rather general class of polynomial ensembles of random matrices  it was first discovered by Johansson 
in (\cite{Jo2})) 
thus contributing a constant term to the r.h.s. of (\ref{predel}), so that
\begin{eqnarray}
\label{clt}
& & E \left( \det(1 + \frac{t^2}{n} \* A^t \* A) \right)^{-1/2}  = E \left ( \exp\left( -\frac{1}{2} \* \sum_{i=1}^n
\ln(1 + t^2\* \lambda_i\*n^{-1}) \right) \right) \nonumber \\
&=& \exp\left( -\frac{1}{2} \* E \big( \sum_{i=1}^n (\ln(1 + t^2\* \lambda_i\*n^{-1}) )\big)  + \frac{1}{8} \* d(t) + o(1) \right),
\end{eqnarray}
where $d(t):=\lim_{n\to \infty} Var\big(\sum_{i=1}^n (\ln(1 + t^2\* \lambda_i\*n^{-1}) )\big)$  and can be explicitely calculated
(see (\cite{Jo2}), Theorem 2.4).

{\bf Remark 5}

Another class of random matrices we are particularly interested in are
Rademacher random matrices
(i.e. square  random matrices with $ \pm 1 $ i.i.d. entries), which we are going to denote by $ \ \mathcal{R} . $
It appears that the questions of the invertibility 
of a Rademacher random matrix and the estimate of the norm of the inverse are of big importance in geometric 
functional analysis (for example in connection with a deterministic construction of Euclidean sections of 
convex bodies). Similar to the previous analysis one can obtain
\begin{equation}
\label{rademacher}
E \ \frac{1}{\sqrt{\det( 1 + t^2 \* \mathcal{R}^t \mathcal{R})}} = 
\frac{1}{(2 \pi)^{n}} \int_{R^{2n}}\prod_{i=1}^{n} du_i dv_i 
 \* \exp\left(- \frac{1}{2} \* \sum_{i=1}^n (u_i^2 + v_i^2) \right) \* \prod_{j,k=1}^{n} 
\cos\Big( t \* u_j\*v_k \Big).
\end{equation}
The fact that a Rademacher matrix $ R $ is invertible with probability going to $ 1 $  as $ n \to \infty $ was 
proved by J.Koml\'{o}s (see e.g. \cite{Bol}, chapter 14). More recently, 
J.Kahn, J.Koml\'{o}s and E.Szemeredi (\cite{KKS}) proved  that the probability that $ R $ is invertible is 
exponentially close to $ 1 .$ To the best of our knowledge there is no known estimate
on the norm of the inverse matrix (which, in our language, corresponds to the estimate  of the smallest 
eigenvalue of $ \ R^tR $).

\section{Complex Matrices with i.i.d. Entries}

In this section we consider the  ensemble of $ n \times n \ $ complex random  matrices $M=A^* A, \ \ 
A= (A_{jk})_{1\leq j \leq m, \ 1 \leq k \leq n} \ $
with the joint distribution of the matrix entries of $A$ given by the formula
\begin{equation}
\label{fyodorov}
\Pr (A) \* dA \* d\overline{A} = \prod_{1 \leq j \leq m, \ 1 \leq k \leq n} \bigl \{ d\Re a_{jk} \* d \Im a_{jk} \* 
\frac{1}{\pi} f(|a_{jk}|^2) \bigr \}. 
\end{equation}

In other words 
$ \{ a_{jk}, \ \ 1 \leq j \leq m, \ 1 \leq k \leq n \} $ 
are independent indentically distributed 
random variables with a distribution depending only 
on the radial component, and $f(x) $ is the 
density of the distribution  of $ |a_{jk}|^2 .$   The ensemble
(\ref{fyodorov}) is a generalization of the standard
Wishart (Laguerre) ensemble which corresponds 
to the choice $ f(x) = e^{-x}.$

In the standard Wishart (Laguerre) case it is known that the smallest eigenvalues are proportional to $ \ \frac{1}{n^2}, \ $
and the (rescaled) $k$-point correlation functions are given in the limit $ \ n\to \infty \ $  by the 
determinants  

\begin{equation}
\label{detcorr}
\rho_k(x_1, \ldots, x_k) = \det(K(x_i, x_j))_{i,j=1, \ldots, k}, \ \ \ k=1,2,3,\ldots, 
\end{equation}

with the Bessel kernel (with $  \ \alpha = 0)$.

\begin{equation}
\label{bessel}
 K^{(\alpha)}(x,y) = \frac{ J_{\alpha}(2 \sqrt{x}) \* \sqrt{y} \* J_{\alpha}'(2 \sqrt{y}) -  
J_{\alpha}(2 \sqrt{y}) \* \sqrt{x} \* J_{\alpha}'(2 \sqrt{x})}{x-y},
\end{equation}
where $J_{\nu}$ is the J-Bessel function, appears as the limit of the rescaled correlation kernel at the hard 
edge in the Laguerre and Jacobi ensembles (see e.g. \cite{For}, \cite{TW3}). 

Ban Arous and P\'{e}ch\'{e}  (\cite{BAP}) following the approach suggested by Johansson (\cite{Jo1}) for Wigner 
matrices have recently shown universality of  the limiting
distribution of the smallest eigenvalues (as well as in the bulk of the spectrum) for a special class of 
sample covariance matrices. 
Their technique requires that  entries of $A$ have a Gaussian component.

One of the possible ways to attack this problem for an
ensemble (\ref{fyodorov})
(assuming that 
all  moments exist, i.e.
and do not grow very fast)
is to study the mathematical expectation of the ratio of 
determinants
\[
Z(\eta_1, \eta_2, \ldots, \eta_k, \mu_1, \ldots, \mu_l) = 
E \frac{ \prod_{i=1}^k \det{(1 + \eta_i^2 \* A^*A)}}{ \prod_{j=1}^l
\det{( 1 + t_j^2 \* A^*A)}} 
\]
for appropriately scaled (large) real numbers
$ \eta_1, \ldots, \eta_k, t_1, \ldots, t_l. \  $ 
For the standard complex Wishart (Laguerre) case these expectation
values were calculated
exactly for any $k,l,n$, see \cite{FA}, \cite{SV} and references
therein, and also used to address
objects interesting in mesoscopic physics \cite{FO}, and Quantum
Chromodynamics\cite{DA}.
 
In particular, one can easily see that
$ \frac{\partial}{\partial \eta} \* Z(\eta, t)|_{\eta=t} = 
E \sum_{i=1}^n \frac{2\eta \* \lambda_i }{t^2 + \lambda_i}, \ \ $ where
$ 0 \leq \lambda_1 \leq \lambda_2 \leq \ldots \leq \lambda_n \ $ are
the eigenvalues of $ A^*A.$ Such an object can be used to extract
the mean eigenvalue density. 
In a similar fashion, by taking partial derivatives of $ Z(\eta_1, \ldots, \eta_k, t_1, \ldots, t_l) $
of higher orders, one can study the correlations of the eigenvalues of $A^*A.$ To show the universality of the 
distribution of the smallest eigenvalues we need to show that local statistical quantities at the edge of the 
spectrum (near the origin) do not depend (in the limit $n\to \infty)$ 
on the second and higher moments of $f.$

The next proposition is analogous to Proposition 1 in the real case.

\begin{proposition}
Let $A$ be a random rectangular $m \times n$ matrix with the probability distribution given by
(\ref{fyodorov}). 
Let $ t_l>0, \ l=1, \ldots, r $ be some positive parameters.
Then the following formula holds
\begin{eqnarray}
\label{complprod}
& & E \left( \prod_{l=1}^r \det(1+ t_i^2\* A^t\*A) \right)^{-1} = \pi^{-r(n+m)/} \ \int_{R^{2r(n+m)}} 
\prod_{l=1}^r d^{n} \Re s^{(l)}\* d^{n} \Im s^{(l)} \*
d^{m} \Re p^{(l)} \* d^{m} \Im p^{(l)}
\nonumber \\
& &   \exp\left( -\sum_{l=1}^r (|s^{(l)}|^2 +|p^{(l)}|^2) \right) 
\prod_{1\leq j \leq n, 1 \leq k \leq m}
G\left(|\sum_{l=1}^r t_l \* s_k^{(l)} \* p_j^{(l)}|^2 \right),
\end{eqnarray}
where $ s^{(l)}=(s^{(l)}_1, \ldots, s^{(l)}_n) $ are $n$-dimensional complex vectors, 
$p^{(l)}=(p^{(l)}_1, \ldots, p^{(l)}_m) $ are $m$-dimensional
complex vectors, $ \ l=1, \ldots, r, $ 
\begin{equation}
\label{G}
G(y)  = \frac{1}{2 \pi} \* \int_{0}^{2 \pi} d \theta
\* \int_0^{+\infty} f(x) \* dx \* \exp \left(2 \* i \* (xy)^{1/2} \* \cos(\theta) \right)=
\int_{0}^{+\infty} dx \* f(x) \* \phi(xy), 
\end{equation}
and
\begin{equation}
\label{phi}
\phi(x)=\sum_{l=0}^{+\infty} \frac{(-1)^l}{(l!)^2} \* x^l = J_0(2 \* x^{1/2}).
\end{equation}
In the special case of single determinant ($r=1$) the formula (\ref{complprod}) can be simplified
\begin{equation}
E \bigl(\det(1 + t^2 \* A^*A)\bigr)^{-1}=
\int_{(0, +\infty)^{m+n}} \prod_{i=1}^n \*  e^{-u_i}\* du_i \* \prod_{j=1}^m e^{-v_i} \* dv_i 
\* \prod_{k,l} G \left( t^2 \* u_k \* v_l \right).
\end{equation}

\end{proposition}

{\bf Remark 6}
As in the section 2 we can consider the case when the  the matrix entries $ \{ |a_{jk}| \} $ are independent but not 
identically distributed
with the densities $ f_{jk}(|x|), \ 1 \leq j \leq m, \ 1 \leq n \leq n. $  The result of Proposition 2 still holds true provided
we replace $G$ in (\ref{complprod}) and (\ref{G}) by $ G_{jk}(y) = \frac{1}{2 \pi} \* \int_{0}^{2 \pi} d \theta
\* \int_0^{+\infty} f_{jk}(x) \* dx \* \exp \left(2 \* i \* (xy)^{1/2} \* \cos(\theta) \right)=
\int_{0}^{+\infty} dx \* f_{jk}(x) \* \phi(xy). $

{\bf Remark 7}

In (\ref{phi})  $ \ J_0(z)= \frac{1}{2 \pi} \* \int_0^{2\pi} \exp(i \* z \* \cos \theta) \* d \theta \ $
is the standard Bessel function (\cite{Wa}).

{\bf Remark 8}

In the special Wishart (Laguerre)
case ( which corresponds to $ f(x)= e^{-x} ) \ $  one has $\ G(y)= e^{-y}. \ $ 
If all moments of $f(x)$ exist and do not grow very fast
one can write
$ G(y) = \sum_{l=0}^{+\infty} \frac{(-1)^l \* \alpha_l }{(l!)^2} \* y^l, $
where
$ \{\alpha_n \}_{n\geq 1} $ are the moments of $f(x).$   

{\bf Proof of Proposition 2.}

Let  $ s= (s_1, \ldots, s_n)^t, \ \ p=(p_1, \ldots, p_m)^t \ $ be complex $n-$ and $m-$ 
dimensional column vectors and $ s^*= (\overline{s_1}, \ldots, \overline{s_n}), \ \ 
p^*= (\overline{p_1}, \ldots, \overline{p_m}).$ 
In what follows $ d^{2n} s $ and $ d^{2m} p $ will stand for $ d^n \Re s \* d^n \Im s $ and
$ d^m \Re p \* d^m \Im p $ correspondingly.
Then
\begin{eqnarray*}
& &  \bigl(\det(1 + t^2 \* A^*A)\bigr)^{-1} = 
\left( \det \left( \begin{array}{cc} 1 & t i \* A\\ t i \* A^* & 1 \end{array} \right) \right)^{-1} \\
& & = \left(\frac{1}{\pi}\right)^{n+m} \*  \int \ d^{2n} s \  d^{2m} p  \ 
\exp \big(-(s^*, p^*) \* B(t) \* (s,p)^t \big)   
 \\
& & = \left(\frac{1}{\pi}\right)^{n+m} \* 
\int \ d^{2n} s \  d^{2m} p \  \exp( - (|s|^2 + |p|^2)) \* \exp \left(
-i \sum_{1 \leq j \leq n, 1 \leq k \leq m} \left (t \* a_{jk} \* p_k \* \overline{s_j} + 
t \* \overline{a_{jk}} \* \overline{p_k} \* s_j \right ) \right),
\end{eqnarray*}
where as before $ B(t)= \left( \begin{array}{cc} Id & it \* A\\ it \* A^* & Id \end{array} \right) \ $ and 
$ \ d^{2n}s \* d^{2m}p = \prod_{i=1}^n d\Re s_i \* d \Im s_i \* \prod_{j=1}^m  \Re p_i \* d \Im p_i. \ $

We can then write down
\begin{eqnarray}
\label{desyat}
& & E \left( \prod_{l=1}^r \det(1+ t_i^2\* A^t\*A) \right)^{-1} = \nonumber \\
& &  \pi^{-r(n+m)} \int_{R^{2r(n+m)}} 
\prod_{l=1}^r d^{2n} s^{(l)}\*
d^{2m} p^{(l)} \ \exp\left( -\sum_{l=1}^r (|s^{(l)}|^2 +|p^{(l)}|^2) \right) \times \nonumber \\ 
& & \prod_{1\leq j \leq n, 1 \leq k \leq m}
E \left(\exp\left(- i \*( a_{kj} \* \sum_{l=1}^r t_l \*s_j^{(l)} \* \overline{p_k^{(l)}} +
\overline{a_{kj}} \* \sum_{l=1}^r t_l \*\overline{s_j^{(l)}} \* p_k^{(l)})
 \right)\right) 
\end{eqnarray}
Let $u= \sum_{l=1}^r t_l \*s_j^{(l)} \* \overline{p_k^{(l)}}. $ Then we can write
\begin{eqnarray}
\label{dven}
& & E \left( \exp \bigl(
- i  (a_{jk} \* u + 
\overline{a_{jk}} \*\overline{u} \bigr) \right)
=  \frac{1}{\pi} \* \ \int  d \Re z \* d \Im z  \* f(|z|^2) \* \exp 
\bigl( - i ( z \* u + \overline{z} \* \overline{u}) \bigr) \nonumber \\
& & = \frac{1}{\pi}  \*  
\int_0^{2\pi} d\theta \* \int_{0}^{+\infty} dr \* r \* f(r^2) \* \exp \bigl(
-  i\*r (e^{i\theta} \* u + e^{-i\theta} \* \overline{u}) \bigr)\nonumber \\ 
& & = \frac{1}{2\pi} \*  
\int_0^{2\pi} d\theta \* \int_{0}^{+\infty} dx \* f(x) \* \exp \bigl(
-  i\* \sqrt{x} (e^{i \* \theta} \* u + e^{-i \* \theta} \* 
\overline{u}) \bigr) \nonumber \\
& & = 
\int_{0}^{+\infty} dx \* f(x) \* \phi(x\*|u|^2),
\end{eqnarray}
where $\phi(x)$ has been defined in (\ref{phi}). Combining (\ref{desyat}) and (\ref{dven}) we arrive at
\begin{eqnarray}
\label{mainform}
& & E \prod_{l=1}^r \bigl(\det(1 + t_l^2 \* A^*A)\bigr)^{-1} = \left(\frac{1}{\pi}\right)^{r(m+n)} \* 
\prod_{l=1}^r \int
\ d^{2n} s^{(l)}  \  d^{2m} p^{(l)}  \ 
\exp \bigl(-  \sum_{j=1}^n |s^{(l)}_j|^2 - \sum_{k=1}^m |p^{(l)}_k|^2) \bigr) \nonumber \times \\
& & \prod_{1\leq j \leq n, 1 \leq k \leq m}
G\left(|\sum_{l=1}^r t_l \* s_j^{(l)} \* p_k^{(l)}|^2 \right).
\end{eqnarray}
In the special case $ r=1$ the formula can be simplified further
\begin{eqnarray}
\label{goodform}
& & E \bigl(\det(1 + t^2 \* A^*A)\bigr)^{-1} = \left(\frac{1}{\pi}\right)^{m+n} \* \int_{R^{2(n+m)}}
\ d^{2n} s  \  d^{2m} p  \
\exp \bigl(-  \sum_{j=1}^n |s_j|^2 - \sum_{k=1}^m |p_k|^2) \bigr) \nonumber \times \\
& & \prod_{1 \leq j \leq n, 1 \leq k \leq m} \ 
\int_0^{\infty} dx \* f(x)  \left\{ \sum_{l=0}^{\infty} \frac{1}{(l!)^2} \* (-x)^l
t^{2l} \* \* |p_k|^{2l} \* |s_j|^{2l} \right\}  \nonumber \\
& & = 2^{m+n} \*  \int_{(0, +\infty)^{m+n}} \prod_{i=1}^n r_i \* dr_i\* \prod_{j=1}^m \rho_j \* d\rho_j 
\* \exp\left(- \sum_{i=1}^n r^2_i - \sum_{j=1}^m \rho^2_j) \right)
\* \prod_{k=1}^n \* \prod_{l=1}^m \*  G \left( t^2 \* r^2_k \* \rho^2_l \right) \nonumber \\
& &  =   \int_{(0, +\infty)^{m+n}} \prod_{i=1}^n \*  e^{-u_i}\* du_i \* \prod_{j=1}^m e^{-v_i} \* dv_i 
\* \prod_{k,l} G \left( t^2 \* u_k \* v_l \right).
\end{eqnarray}

In the Wishart case  one can simplify things even further. The calculations are very similar to the real 
Wishart case considered in Remark 4 (section 2). Since typical eigenvalues of $A^*\*A$ are of the order of $n$
the scaling of $t^2$ by a factor $1/n$ allows us to study the limiting distribution of the eigenvalues. It 
follows from (\ref{goodform}) that
we are left  with the task of evaluating the integral
$ E \bigl(\det(1 + \frac{t^2}{n} \* A^*A)\bigr)^{-1} 
\int_{(0, +\infty)^{2n}} \prod_{i=1}^n  e^{-u_i}\* du_i\*  \prod_{j=1}^m e^{-v_j} \* dv_j 
\* \exp \left( -\frac{t^2}{n} \* \sum_{k=1}^n u_k \* \sum_{l=1}^m v_l \right),  $ which can be reduced to

\begin{eqnarray}
\label{kto}
& & \int_{(0, +\infty)^{m+n}} \prod_{i=1}^n  e^{-u_i}\* du_i\* \prod_{j=1}^m e^{-v_j} \* dv_j
\* \exp \left( -\frac{t^2}{n} \* \sum_{k=1}^n u_k \* \sum_{l=1}^m v_l \right)  \nonumber \\
& & = \* \int_{(0, +\infty)^{n}} \prod_{i=1}^n  e^{-u_i} \* du_i \* \left( 1 + \frac{t^2}{n} \* 
\sum_{l=1}^n u_l \right)^{-m} \nonumber \\
& & =\int_0^{+\infty} \frac{z^{n-1}}{(n-1)!} \* e^{-z} \* \left(1+ \frac{t^2}{n} 
\* z \right)^{-m} \* dz \nonumber \\
& & = \frac{n^n}{\Gamma(n)} \* \int_0^{+\infty} z^{n-1} \* e^{-nz} \* 
\left(1+  t^2 \* z \right)^{-m} \ dz = 
\frac{n^n}{\Gamma(n)}\int_0^{+\infty} e^{-n\*\mathcal{L}(z)} \* z^{-1} \* dz
\end{eqnarray}
where $ \ \mathcal{L}(z)=  z + \frac{m}{n} \*\ln(1 +t^2\*z)
-\ln z. \ $
It is remarkable that the formulas in the complex case 
are identical to those in the real case (\ref{gaussnol}-\ref{gaussodin}) modulo trivial
change of parameters. We thus proved the following result.

\begin{lemma}
Let $m$ and $n$ be positive integers and $z \in \Bbb C \setminus (-\infty,0)$. Then
\begin{equation}
\label{equality}
E_{2m,2n, real} \left(\det(1+\frac{z}{2} \*A^t\*A)\right)^{-1/2} = E_{m,n, complex} 
\left(\det(1+z\*A^*\*A)\right)^{-1}
\end{equation}
where at the l.h.s. we have the mathematical expectation with respect to the ensemble of rectangular
$2m\times 2n$ real matrices $A$ with  i.i.d. standard Gaussian entries (standard real Wishart ensemble), and the
r.h.s. we have the mathematical expectation with respect to the ensemble of rectangular
$m\times n$ complex matrices $A$ with i.i.d. standard Gaussian entries 
(standard complex Wishart ensemble)
\end{lemma}

As before, it was enough to prove the result for positive real $z=t^2. $
We remind the reader that in the standard real Wishart case all entries $ \{ a_{j,k} \} $ are i.i.d. N(0,1) random variables,
and in the standard complex Wishart case all entries $ \{ \Re a_{j,k}, \Im a_{j,k} \}$ are i.i.d. N(0, 1/2) random variables
(so in both cases $ E |a_{j,k}|^2 =1 $).

{\bf Remark 9}

If all moments of $f(x)$ exist then under some technical conditions
the asymptotics of $ \ E \left(\det\left(1 + t^2 \* A^*A \right) \right)^{-1} \ $ 
in the global regime 
depends
on the first moment and second moments of $ \ f(x)dx, \ $  (i.e.  on the second and fourth moments
of the matrix entries $ A_{kl}. \  $  This phenomena is known in random matrix 
theory: for example in the case of Wigner random matrices, the limiting distribution of a global linear 
statistics $ \ Tr \ h(A) \ - E (Tr \ h(A)) \ $ (where $h$ is a test function, say a polynomial, and $ A $ is a 
random Wigner matrix normalized so that a typical eigenvalue is of the order of a constant) depends on the 
second and 
fourth moments of the matrix entries (see e.g. \cite{KKP}, \cite{SiS1}). 
It is conjectured (and in a few interesting special cases verified) that
in the local regime the dependence on the fourth
moment  goes away.

\section{Appendix}

A Poisson random point process on the positive half-axis with the locally integrable intensity function $ \rho(x) $ is defined in such a way
that the counting functions (e.g. numbers of particles) in the disjoint intervals $ I_1, \ldots, I_k $ are independent Poison random 
variables with the parameters $ \int_{I_j} \rho(x) \* dx, \ \ j=1, \ldots, k.$ Equavalently, one can define the Poisson random point by 
requiring that the $k$-point correlations functions are given by the products of one-point correlation functions (intensities), i.e.
$\rho_k(x_1, \ldots, x_k)= \prod_{j=1}^k \rho(x_j).$ 

Let $f: (0, +\infty) \to \Bbb C$ be a test function with a nice behavior at the origin and infinity. Then
\begin{eqnarray}
\label{poisson}
& & {\bf E} \prod_{i=1}^{\infty} (1+f(x_i)) = 1 +\sum_{k=1}^{\infty} {\bf E} \sum_{1 \leq i_1 < i_2 <\ldots < i_k} 
\prod_{j=1}^k f(x_{i_j}) = 
\sum_{k=0}^{\infty} \frac{1}{k!} \int_{(0, +\infty)^k} \prod_{j=1}^k f(x_j) \* \rho_k(x_1, \ldots, x_k) dx_1 \cdots dx_k \nonumber \\
&=&
\sum_{k=0}^{\infty} \frac{1}{k!} \left( \int_{(0, +\infty)} f(x)\* \rho(x) dx \right)^k=
\exp\left( \int_{(0, +\infty)} f(x)\* \rho(x) dx \right)
\end{eqnarray}
If the test function $f(x) $ equals $ (1 +z\*x)^{-1/2} -1, $ and $ \rho(x)= \frac{1}{\pi\* x^{3/2}} $ we have
$$\int_{(0, +\infty)} f(x)\* \rho(x) dx = \int_{(0, +\infty)} ((1 +z\*x)^{-1/2} -1) \* \frac{1}{\pi\* x^{3/2}} dx = 
-\frac{2}{\pi} \* \sqrt{z}, $$  which is exactly the exponent in (\ref{glavforma}). For random Schr\"{o}dinger operators the Poisson 
statistics of the eigenvalues in the localization regime was first proved by Molchanov in \cite{Mo} (see also \cite{Mi}).

\noindent

{\bf Acknowledgements.}\, It is a pleasure to thank Boris Khoruzhenko for very useful comments.

\def\cmp{{\it Commun. Math. Phys.} }

\end{document}